\documentclass[ukrainian,12pt]{article}
\usepackage{amsmath}
\usepackage{amsbsy}
\usepackage{amssymb}
\usepackage{amscd}
\usepackage{a4wide}
\usepackage[T2A]{fontenc}
\usepackage[cp1251]{inputenc}
\usepackage[ukrainian,english]{babel}
\usepackage{epic}
\newenvironment{proof}{\par\vspace{1ex}\par {\em Proof:}\hspace{0.5em}}{\noindent$\bigtriangleup$\par\vspace{1ex}\par }

\newtheorem{theorem}{Theorem}

\newtheorem{definition}{Definition}

\DeclareMathOperator{\rad}{Rad} 
 \DeclareMathOperator{\rep}{Rep}
 \DeclareMathOperator{\image}{\Im m}
 \DeclareMathOperator{\ob}{Ob}

\begin{document}

\title{ The spectral  problem and algebras
associated with extended Dynkin graphs. I. }
\author{\footnotesize STANISLAV KRUGLJAK, STANISLAV POPOVYCH, YURII SAMOILENKO \footnote{Insitute of Mathematics of National Academy of Sciences, Kyiv, Ukraine}}

\maketitle

\newcommand{\ba}[1]{\ensuremath{e_{\phi(#1)}}}
\newcommand{\sca}[2]{\ensuremath{\langle #1,#2 \rangle}}
\newcommand{\algebra}[1]{\ensuremath{\mathbb{C}\langle #1 \rangle}}
\newcommand{\algp}[1]{\ensuremath{\mathcal{P}_{n,#1}}}
\newcommand{\pabo}[1]{\ensuremath{\mathcal{P}_{n,abo,#1}}}
\newcommand{\g}[2]{\ensuremath{\gamma_{#1}^{(#2)}}}
\newcommand{\ma}[2]{\ensuremath{\alpha_{#1}^{(#2)}}}
\numberwithin{equation}{section}
 \abstract{  There is a  connection between *-representations of algebras
associated with graphs and the problem about the  spectrum of a
sum of  Hermitian operators (spectral problem). For algebras
associated with extended Dynkin graphs we give an explicit
description of the parameters for which   there are
$*$-representations and an algorithm for constructing these
representations.


 \medskip\par\noindent
AMS SUBJECT CLASSIFICATION: 16W10, 16G20, 47L30
  }
\vspace{10mm}

\noindent {\Large \bf Introduction.} \vspace{5mm}

 {
 1. Let  $A_1$, $A_2$, $\ldots$, $A_n$  be Hermitian  $m\times m$
 matrices with given
eigenvalues: $\tau(A_j)= \{\lambda_1(A_j)\ge \lambda_2(A_j)\ge
\ldots \ge \lambda_m(A_j) \}$. The well-known classical problem
about the spectrum of a sum of two Hermitian matrices (Horn's
problem) is to describe a connection between
 $\tau(A_1), \tau(A_2), \tau(A_3)$  such that $A_1+A_2=A_3$. In
 more symmetric setting one can seek for a connection between $\tau(A_1), \tau(A_2), \ldots,
 \tau(A_n)$ necessary and sufficient for the existence of Hermitian
 operators $A_1+A_2+\ldots+A_n=\gamma I$ for a fixed $\gamma\in
 \mathbb{R}$.

A recent solution of this problem  (see ~\cite{fulton, klychko}
 and others) gives a complete description of
possible $\tau(A_1), \tau(A_2), \ldots,
 \tau(A_n)$ in terms of linear
inequalities.  Note that the number of necessary inequalities
increases with $m$.

\noindent 2. In \cite{KPS1, KPS} was considered the following
modifications of the problem mentioned above,  called henceforth
the {\it spectral problem} (resp. {\it strict spectral  problem}).
We will consider bounded linear Hermitian operators $A_1$, $A_2$,
$\ldots$, $A_n$ on a separable Hilbert space $H$. For an operator
$X$ denote by $\sigma(X)$ its spectrum. Let $M_1, M_2, \ldots,
M_n$ be given closed subsets of $\mathbb{R}$ and
$\gamma\in\mathbb{R}$. The problem consists of the following: 1)
to determine whether there are Hermitian operators $A_1$, $A_2$,
$\ldots$, $A_n$ on $H$ such that $\sigma(A_j)\subseteq M_j$
(respectively $\sigma(A_j) = M_j$) ($1\le j\le n$), and
$A_1+A_2+\ldots+A_n=\gamma I$? 2) if the answer is in the
affirmative, to give a description (up to unitary equivalence) of
the operators. In this work the sets $M_1, M_2, \ldots, M_n$ will
be finite. Note that even for finite $M_k$ the second part of the
problem can be very complicated if $|M_k|$ is large enough.

The essential difference with Horn's classical  problem is that we
do not fix the dimension of  $H$   (it may be finite or infinite)
  and the spectral multiplicities. It seems that the  solution of the spectral
   and strict spectral  problems could not be deduced directly from the Horn inequalities,
   since the number of necessary inequalities increases with $m$.

\noindent 3. Spectral problem could be reformulated in terms of
*-representations of *-algebras.  Namely, let
$\alpha^{(j)}=(\alpha_1^{(j)},\alpha_2^{(j)},\ldots,
\alpha_{m_i}^{(j)})$ ($1\le j\le n$)  be vectors with positive
strictly decreasing coefficients.  Put  $M_j=\alpha^{(j)}$. Let us
consider the associative algebra defined by the following
generators and relations (see.~\cite{popsam}):
\begin{gather*}
\mathcal{A}_{M_1,  \ldots, M_n, \gamma}=\mathbb{C}\langle
p_1^{(1)}, p_2^{(1)}, \ldots, p_{m_1}^{(1)}, p_1^{(2)},
p_2^{(2)},\ldots, p_{m_2}^{(2)}, \ldots, p_1^{(n)},
p_2^{(n)},\ldots,
p_{m_n}^{(n)} |\\
 p_k^{(i)2}=p_k^{(i)}, \sum_{i=1}^{n} \sum_{k=1}^{m_i}\alpha_k^{(i)}
  p_k^{(i)} = \gamma e, \\ \sum_{k=1}^{m_i}p_k^{(i)} =e,
  p_j^{(i)}p_k^{(i)} = 0, j,k=1,\ldots, m_i, j\not=k, i=1,\ldots,
  n
  \rangle.
\end{gather*}
 Here $e$ is the identity of the algebra.
 This is a  $*$-algebra, if we declare all generators to be
self-adjoint.

A $*$-representation $\pi$ of $\mathcal{A}_{M_1,  \ldots, M_n,
\gamma}$  determines  an $n$-tuple  of non-negative operators
$A^{(j)}=\sum_{k=1}^{m_i} \alpha_k^{(i)} P_k^{(i)}$,
  where each of the families of orthoprojections, $\{
P_i=\pi(p_i)$, $i=1,\ldots, k  \}$ forms a resolution of the
identity and  such that  $A^{(1)}+\ldots+ A^{(n)}=\gamma I$.
Moreover, $\sigma(A^{(j)})\subseteq M_j$. And viceversa, any
$n$-tuple of Hermitian matrices $A^{(1)}$, $A^{(2)},\ldots$,
$A^{(n)}$ such that $A^{(1)}+\ldots+ A^{(n)}=\gamma I$ and
$\sigma(A^{(j)})\subseteq M_j$ determines a representation $\pi$
of $\mathcal{A}_{M_1,  \ldots, M_n, \gamma}$ by taking
$\pi(p_k^{(j)})$ to be the spectral projections corresponding to
eigenvalues $\alpha_k^{(j)}$ respectively.
 So in terms of $*$-representations, the
spectral problem is a problem consisting of the following two
parts: 1) a description of the set $\Sigma_{m_1,m_2,\ldots,m_n}$
of the parameters $\alpha_k^{(j)}$, $\gamma$ for which there exist
$*$-representations of $\mathcal{A}_{M_1,  \ldots, M_n, \gamma}$.
   2) a description of
*-representations $\pi$ of the *-algebra $\mathcal{A}_{M_1,
M_2,\ldots, M_n, \gamma}$.

 \noindent 4. A natural way to try to  solve the spectral  problem is
 to describe all irreducible $*$-representations up to unitary
 equivalence and then all *-representations as sums or
 direct integrals of irreducible representations.
Obviously if there is  a  *-representation of the algebra
$\mathcal{A}_{M_1,  \ldots, M_n \gamma}$ then there is its
irreducible *-representation. Hence the set
$\Sigma_{m_1,m_2,\ldots,m_n}$ coincides with the set of parameters
for which algebra $\mathcal{A}_{M_1,  \ldots, M_n, \gamma}$ has at
least one irreducible *-representation;
  The second
 part of the spectral problem could be formulated in the following
 way: find the formulae for the irreducible *-representations
 of the algebra $\mathcal{A}_{M_1,  \ldots, M_n, \gamma}$ for
 parameters $(\lambda_k^{(j)},\gamma)\in \Sigma_{m_1,\ldots, m_n}$
 or at least present an algorithm to construct such
 representations. The strict spectral problem could also be
 reformulated in terms of representation theory of these algebras
 but we  will not discuss it here.

 \noindent  5. An key step in solving  spectral problem is
   to describe the irreducible
non-degenerate representations.  Let us call a  $*$-representation
$\pi$  of the algebra
 $\mathcal{A}_{M_1,  \ldots, M_n, \gamma}$  {\it non-degenerate}
if  $\pi(p_k^{(j)})\not=0$ for all $k$, $j$.
 Consider the following set:
    $\Sigma_{m_1,\ldots,m_n}^{n.-d.}=
  \{ (\{\lambda_k^{(j)}\},\gamma)|$ { \it there is a non-degenerate *-representation of }
  $\mathcal{A}_{M_1,  \ldots, M_n, \gamma}\}$;
  which depends only on   $(m_1,\ldots,m_n)$.
 Every irreducible representation of the algebra $\mathcal{A}_{M_1,
M_2, \ldots, M_n, \gamma}$   irreducible non-degenerate
*-representation of an algebra $\mathcal{A}_{\widetilde{M}_1,
\ldots, \widetilde{M}_n, \gamma}$ for some subsets
$\widetilde{M}_j\subset M_j$. Hence $(M_1,  \ldots, M_n,
\gamma)\in \Sigma_{m_1,\ldots, m_n}$ if there exist
$(\widetilde{M}_1, \ldots, \widetilde{M}_n,\gamma)\in
\Sigma_{|\widetilde{M}_1|,\ldots, |\widetilde{M}_n|}^{n.-d.}$.
Thus the description of $\Sigma_{m_1,\ldots, m_n}$ follows from
the description of $\Sigma_{k_1,\ldots, k_n}^{n.-d.}$ where
$k_j\le m_j, 1\le j\le n$.

\noindent 6.   With an  integer vector    $(m_1,\ldots,m_n)$ we
will
  associate a non-oriented  star-shape  graph  $G$
  with $n$ branches of the lengths   $m_1$, $m_2$, $\ldots$,
  $m_n$  stemming from a  single root. The graph $G$ and vector
  $\chi= (\alpha_1^{(1)}, \alpha_2^{(1)}, \ldots, \alpha_{m_1}^{(1)},
  \alpha_{1}^{(2)}, \ldots, \alpha_{m_2}^{(2)}, \ldots,
  \alpha_{1}^{(n)}, \ldots, \alpha_{m_n}^{(n)}, \gamma)$ completely
  determine the algebra $\mathcal{A}_{M_1,  \ldots, M_n,
  \gamma}$   so we will use the following notation $ \mathcal{A}_{G,
  \chi}$.

  Henceforth we will denote $\Sigma$ by  $\Sigma(G)$ where $G$ is the
  tree mentioned above.
The  spectral problem   for operators on a Hilbert space can be
reformulated in the following way: 1) for a given graph  $G$
describe the set  $\Sigma(G)$; 2) describe non-degenerate
representations $\mathcal{A}_{M_1,  \ldots, M_n, \gamma}$ up to
unitary equivalence.

If the  graph is a Dynkin graph(extended Dynkin graph), the
problem is greatly simplified. The algebras $\mathcal{A}_{M_1,
\ldots, M_n, \gamma}$ associated with Dynkin graphs  (resp.
extended Dynkin graphs) have a more simple structure than in other
cases. In particular, the algebras $\mathcal{A}_{M_1,  \ldots,
M_n, \gamma}$ are finite dimensional (resp. have polynomial
growth) if and only if the associated graph is a Dynkin graph
(resp. an extended Dynkin graph) (see~\cite{melit}). As shown
in~\cite{roiter} irreducible representations of the algebras
associated with Dynkin graphs exist only in certain dimensions
that are bounded from above. In \cite{KPS1, KPS, Zav} we have
given a complete description of $\Sigma(G)$ for all Dynkin graphs
$G$ and an algorithm for finding  all irreducible representations.
In present paper we do the same for extended Dynkin graphs.

\noindent 7. To solve the spectral problem for extended Dynkin
graphs we will follow the following scheme:

1.) As it was mentioned above it in suffices to describe the sets
$(M_1, \ldots, M_n,\gamma)$ for which there exist irreducible
non-degenerate representations of the algebra $\mathcal{A}_{M_1,
\ldots, M_n, \gamma}$. Any such representation is finite
dimensional (see~\cite{ostrovskij}).

2.) Using the connection between representation theory of
$\mathcal{A}_{M_1, \ldots, M_n, \gamma}$ and locally-scalar
representations of the associated graph (see s.~\ref{repconnec})
we can find generalized dimensions of such representations of the
algebra since they correspond to the positive roots of the
corresponding root system. The root can be real or imaginary.

3.) In the forme case such representation can be obtained using
Coxeter functors (see s.~\ref{s1}) starting from the simplest
ones. The corresponding sets of parameters in case of graph
$\widetilde{E}_6$ are described in theorems~\ref{ineqe6},
\ref{ineqe6_2} and \ref{ineqe6_3}.

4.) In the latter case the parameters of the algebra belong to a
certain hyperplane (see s.\ref{repconnec}). In this case the
dimension of irreducible representation is the minimal imaginary
root (or real root which is considered before) which is unique for
each extended Dynkin graph. Since in this case the dimension is
fixed the solution of the spectral problem could be obtained by
direct application of Horn's inequalities. For example for
$\widetilde{E}_6$ graph corresponding inequalities are witten out
in theorem~\ref{E6_hyp}.

Since the description of   $*$-representations of the $*$-algebras
associated with extended Dynkin graphs   can be reduced to a
description of non-degenerate irreducible representations, in this
article we essentially solve the spectral  problem for the
algebras associated with extended Dynkin graphs.

In the present article we develop the general setup for solving
the spectral problem for extended Dynkin graphs illustration the
concepts on example of $\widetilde{E}_6$. In further publications
we will present explicit solution for all extended Dynkin grapsh.

\section{
Locally-scalar graph representations and representations of the
algebras generated by orthoprojections.}\label{s1}

The main tool for our classification is Coxeter functors for
locally-scalar graph representations. They allow one to construct
all representations starting from the simplest ones which
correspond to the vertices of the graph. First we will recall a
connection between category of *-representation of algebra
$\mathcal{A}_{M_1, \ldots, M_n, \gamma}$ associated with the graph
$G$ and locally-scalar representations of the  graph  $G$. For
more details see~\cite{KPS1}. Note that Coxeter functors could be
constructed directly for the categories of representations of
algebras $\mathcal{A}_{M_1, \ldots, M_n, \gamma}$
(see~\cite{popsam}) but the simplest representations of a graph
does not correspond to a representation of corresponding algebra.
This force us to use graph representation terminology and
techniques.

 Henceforth we will use definitions, notations and
results
 about representations of graphs in the category of
Hilbert spaces found in ~\cite{roiter}. Let us recall some of
them.

A graph $G$ consists of a set of vertices $G_v$ a set of edges
$G_e$ and a map $\varepsilon$ from $G_e$ into the set of one- and
two-element subsets of $G_v$ (the edge is mapped into the set of
incident vertices). Henceforth  we consider connected finite
graphs without cycles (trees). Fix a decomposition of $G_v$
 of the form   $G_v= {\overset{\circ}{G}}_v \sqcup
{\overset{\bullet}{G}}_v$ (unique up to permutation) such that for
each  $\alpha\in G_e$ one of the vertices from
$\varepsilon({\alpha})$ belongs to ${\overset{\circ}{G}}_v$ and
the other to ${\overset{\bullet}{G}}_v$.  Vertices in
${\overset{\circ}{G}}_v$ will be called even, and those in the set
${\overset{\bullet}{G}}_v$ odd.
 Let us recall the definition of a  representation  $\Pi$ of a  graph $G$
 in the category of Hilbert spaces  $\mathcal{H}$. Let us associate with each vertex  $g\in G_v$  a Hilbert space
   $\Pi(g)= H_g\in \text{Ob}
\mathcal{H}$, and with each edge  $\gamma\in G_e$ such that
 $\varepsilon(\gamma) =\{g_1,g_2\}$  a pair of mutually adjoint
operators
  $\Pi(\gamma)=\{ \Gamma_{g_1,g_2}$, $\Gamma_{g_2,g_1}\}$, where
  $\Gamma_{g_1,g_2}:H_{g_2}\to H_{g_1}$. We now construct a category
  $\rep (G,\mathcal{H})$.  Its objects are the representations of the
graph  $G$  in  $\mathcal{H}$. A  morphism  $C: \Pi\to
\widetilde{\Pi}$ is a family   $\{C_g\}_{g\in G_v}$ of operators
  $C_g: \Pi(g)\to \widetilde{\Pi}(g)$ such that the following diagrams
  commute  for all edges   $\gamma_{g_2,g_1}\in G_e$:
\begin{equation*}
\begin{CD}
H_{g_1} @>\Gamma_{g_2,g_1}>>H_{g_2} \\  @V C_{g_1} VV  @VV
C_{g_2}V
\\ \widetilde{H}_{g_1}
@>\widetilde{\Gamma}_{g_2,g_1}>>\widetilde{H}_{g_2}
\end{CD}
\end{equation*}

 Let $M_g$ be the set of vertices connected with  $g$ by an edge.
 Let us define the  operators
\[
A_g= \sum_{g'\in M_g} \Gamma_{gg'}\Gamma_{g'g}.
\]

A representation    $\Pi$ in  $\rep (G,\mathcal{H})$  will be
called {\it locally-scalar} if all operators $A_g$ are scalar,
$A_g= \alpha_g I_{H_g}$. The full subcategory $\rep
(G,\mathcal{H})$,
  the objects  of which are locally-scalar representations, will
be denoted by $\rep G$ and  called the category of locally-scalar
representations of the   graph $G$.

Let us denote by  $V_{G}$ the real vector space consisting of sets
$x=(x_g)$ of real numbers $x_g$, $g\in G_v$. Elements $x$ of $V_G$
we will call  $G$-vectors. A vector $x=(x_g)$ is called positive,
$x>0$, if $x\not=0$ and  $x_g\ge 0$ for all $g\in G_v$. Denote
$V_G^+= \{x\in V_G| x>0\}$.  If $\Pi$ is a finite dimensional
representation of the graph $G$ then the  $G$-vector $(d(g))$,
where $d(g)= \dim \Pi(g)$ is called the {\it dimension} of $\Pi$.
If $A_g=f(g) I_{H_g}$ then the $G$-vector $f=(f(g))$ is called the
{\it character } of the locally-scalar representation   $\Pi$ and
$\Pi$ is called the $f$-representation in this case. The {\it
support} $G_v^{\Pi}$ of
 $\Pi$ is $\{ g\in G_v| \Pi(g) \not= 0 \}$. A representation   $\Pi$
is {\it faithful} if $G_v^{\Pi}= G_v$. A character of the
locally-scalar representation $\Pi$ is uniquely defined on the
support  $G_v^{\Pi}$ and non-uniquely on its complement. In the
general case, denote by $\{f_\Pi \}$ the set of characters of
$\Pi$. For each vertex $g\in G_v$, denote by  $\sigma_g$ the
linear operator on  $V_G$ given by the formulae:
 \begin{gather*}
 (\sigma_g x)_{g'} = x_{g'}\ \text{if}\ g'\not=g,\\
(\sigma_g x)_{g} = -x_{g} +\sum_{g'\in M_g}x_{g'}.
 \end{gather*}
The mapping   $\sigma_g$  is called the {\it reflection} at the
vertex $g$. The composition of all reflections at odd vertices is
denoted by   $ \overset{\bullet}{c}$ (it does not depend on the
order of the factors), and at all even vertices by $
\overset{\circ}{c}$. A Coxeter transformation is
$c=\overset{\circ}{c}\overset{\bullet}{c}$, $c^{-1}=
\overset{\bullet}{c}\overset{\circ}{c}$.  The transformation
$\overset{\bullet}{c}$ ($\overset{\circ}{c}$)  is called an  odd
(even) Coxeter map.
  Let us adopt the  following notations for compositions of the  Coxeter maps:
   $\overset{\bullet}{c}_k= \ldots
\overset{\bullet}{c}\overset{\circ}{c}\overset{\bullet}{c}$ ($k$
factors), $\overset{\circ}{c}_k= \ldots
\overset{\circ}{c}\overset{\bullet}{c}\overset{\circ}{c}$ ($k$
factors), $k\in\mathbb{N}$.

Any real function  $f$ on $G_v$ can be identified  with a
$G$-vector $f=(f(g))_{g\in G_v}$. If $d(g)$  is the dimension of a
locally-scalar graph representation  $\Pi$, then
\begin{gather}
\overset{\circ}{c}(d)(g)=\begin{cases} -d(g) +  \sum_{g'\in
M_g}d(g'), &\text{if}\  g\in \overset{\circ}{G}_v,\\
d(g), & \text{if}\ g\in \overset{\bullet}{G}_v,
\end{cases}\\
\overset{\bullet}{c}(d)(g)=\begin{cases} -d(g) +  \sum_{g'\in
M_g}d(g'), &\text{if}\  g\in \overset{\bullet}{G}_v,\\
d(g), & \text{if}\ g\in \overset{\circ}{G}_v.
\end{cases}
\end{gather}

 For $ d\in Z_G^+$ and  $f\in V_G^+$, consider the full subcategory
  $\rep (G,d,f)$ in $\rep
 G$ (here $Z_G^+$ is the set of positive integer  $G$-vectors),
with the set of objects  $Ob \rep (G, d, f)= \{ \Pi| \dim \Pi (g)
=d(g), f\in\{f_\Pi\}\}$. All representations
    $\Pi$ from $\rep (G,d,f)$ have the same support  $X=X_d=G_v^\Pi=\{g\in G_v|
  d(g)\not=0\}$.   We will consider these categories only if   $(d,f)\in
  S=\{ (d,f)\in Z_G^+\times V_G^+ | d(g) + f(g) >0, g\in G_v \}$.
   Let  $\overset{\circ}{X}= X\cap \overset{\circ}{G}_v$,
    $\overset{\bullet}{X}= X\cap \overset{\bullet}{G}_v$.  $\rep_{\circ}(G, d, f) \subset \rep(G, d,
    f)$  ( $\rep_{\bullet}(G, d, f) \subset \rep(G, d,
    f)$) is the full subcategory with objects  $(\Pi, f)$ where
      $f(g)>0$  if  $g\in \overset{\circ}{X}$ ($f(g)>0$  if  $g\in
    \overset{\bullet}{X}$). Let $S_0=\{(d,f)\in S|f(g)>0 \text{ if }
g\in \overset{\circ}{X}_d\}$, $S_\bullet=\{(d,f)\in S|f(g)>0
\text{ if } g\in \overset{\bullet}{X}_d\}$

Put
\begin{gather}
\overset{\bullet}{c}_d(f)(g)=\overset{\circ}{f}_d(g)=
\begin{cases}
 \overset{\bullet}{c}(f)(g), &\text{if}\  g\in \overset{\bullet}{X}_d,\\
f(g), & \text{if}\ g\not\in \overset{\bullet}{X}_d,
\end{cases}\\
\overset{\circ}{c}_d(f)(g)=\overset{\bullet}{f}_d(g)=
\begin{cases}
 \overset{\circ}{c}(f)(g), &\text{if}\  g\in \overset{\circ}{X}_d,\\
f(g), & \text{if}\ g\not\in \overset{\circ}{X}_d.
\end{cases}.
\end{gather}

 Let us denote  $\overset{\bullet}{c}_d^{(k)}(f)=
 \ldots \overset{\bullet}{c}_{\overset{\circ}{c}_2(d)}
 \overset{\circ}{c}_{\overset{\circ}{c}(d)}
 \overset{\bullet}{c}_d(f)$ ($k$  factors),
 $\overset{\circ}{c}_d^{(k)}(f)=
 \ldots \overset{\circ}{c}_{\overset{\bullet}{c}_2(d)}
  \overset{\bullet}{c}_{\overset{\bullet}{c}(d)}
  \overset{\circ}{c}_d(f)$ ($k$  factors).
The even and odd Coxeter reflection functors are defined
in~\cite{roiter},
  $\overset{\circ}{F}: \rep_{\circ} (G, d, f) \to
  \rep_{\circ} (G, \overset{\circ}{c}(d), \overset{\circ}{f}_d)$
   if  $(d,f)\in S_\circ$,
$\overset{\bullet}{F}: \rep_{\bullet} (G, d, f) \to
  \rep_{\bullet} (G, \overset{\bullet}{c}(d), \overset{\bullet}{f}_d)$
   if  $(d,f)\in S_\bullet$; they are equivalences of the  categories.
Let us denote
   $\overset{\circ}{F}_k(\Pi)=
 \ldots \overset{\circ}{F}\overset{\bullet}{F}\overset{\circ}{F}(\Pi)$ ($k$  factors),
 $\overset{\bullet}{F}_k(\Pi)=
 \ldots \overset{\bullet}{F}\overset{\circ}{F}\overset{\bullet}{F}(\Pi)$ ($k$
 factors),  if the compositions exist. Using these functors, an analog of Gabriel's theorem
for graphs and their locally-scalar representations has been
proven in~\cite{roiter}. In particular, it has been proved that
any locally-scalar graph representation decomposes into a direct
sum (finite or infinite) of finite dimensional indecomposable
representations, and all indecomposable representations can be
obtained by odd and even Coxeter reflection functors starting from
the simplest representations  $\Pi_g$ of the graph  $G$
 ($\Pi_g(g)=\mathbb{C}, \Pi_g(g')= 0$  if $g\not=g';\ g,g'\in G_v
 $).

\section{Root systems associated with extended Dynkin diagrams.}
 Let us recall a few facts about root systems.
 Let $G$ be a simple connected graph. Then its {\it Tits form}
 $$q(\alpha)= \sum_{i\in G_v} \alpha_i^2 -\frac{1}{2} \sum_{\beta\in G_e,
   \{i,j\}=\epsilon(\beta)} \alpha_i\alpha_j  (\alpha \in V_G).
 $$
The {\it symmetric biliniar form} $(\alpha, \beta) = q(\alpha +
\beta ) - q(\alpha) -q(\beta)$. Vector $\alpha \in V_G$ is called
 sincere if each component is non-zero.

 It is well known that for Dynkin graphs (and only for them)
 bilinear form $(\cdot, \cdot)$ is positive definite. The form is
 positive semi-definite for extended Dynkin graphs. And in the
 letter case $\rad q= \{v | q(v)=0 \}$ is equal to $\mathbb{Z}
 \delta$ where  $\delta$  is a minimal imaginary root. For other
 graphs (which are neither Dynkin nor extended Dynkin) there are
 vectors $\alpha \ge 0$ such that $q(\alpha)<0$ and $(\alpha, \epsilon_j)\le
 0$ for all $j$.

 If $G$ is an extended Dynkin graph a vertex $j$ is called
 {\it extending } if $\delta_j =1$. The graph obtained by deleting
 extending vertex is the corresponding Dynkin graph. The set of
 {\it roots} is $\Delta = \{ \alpha \in V_G | \alpha_i \in \mathbb{Z} \text{ for all } i\in G_v, \alpha\not=0, q(\alpha)\le 0
 \}$. A root $\alpha$ is {\it real} if $q(\alpha) = 1$ and
 {\it imaginary} if $q(\alpha) =0$. Every root is either positive
 or negative. For our classification purposes we will need the
 following fact (see~\cite{cw}): for an extended Dynkin graph the
 set $\Delta\cup \{0\}\slash \mathbb{Z} \delta$ is finite. Moreover, if $e$ is an
 extending vertex then the set $\Delta_f =\{\alpha \in \Delta\cup\{0\} |
 \alpha_e=0 \}$ is a complete set of representatives of the cosets
 from $\Delta\cup \{0\} \slash \mathbb{Z} \delta$.
 If $\alpha$ is a root then $\alpha+\delta$ is again a root. We
 will call the coset $\alpha +\delta\mathbb{Z}$ a $\delta$
 -series. If $\alpha$ is a root then its images under the action
 of the group generated by $\overset{\circ}{c}$ and
 $\overset{\bullet}{c}$ will be called a Coxeter series or
 $C$-series for short. It turns out that each $C$-series
 decomposes into a finite number of $\delta$-series of roots. We
 will use this decompositions to give an explicit formulae for
 generalized dimensions and corresponding parameters of algebras.

Let  $G$ be the extended Dynkin graph $\widetilde{E_6}$,
\setlength{\unitlength}{.7mm}
\begin{picture}(62,25)(-1,-1)
\thicklines
 \drawline(19,0)(19,9)
 \multiputlist(0,0)(10,0){\circle*{2},\circle{2},\circle*{2},\circle{2},\circle*{2}}
 \multiputlist(20,10)(0,10){\circle{2},\circle*{2}}
 \drawline(0,0)(8,0)
\drawline(10,0)(19,0)
\drawline(20,0)(28,0)\drawline(30,0)(38,0)\drawline(19,11)(19,19)
\put(0,1){$g_1$} \put(10,1){$g_2$} \put(20,1){$g_0$}
\put(30,1){$g_4$} \put(40,1){$g_3$} \put(21,10){$g_6$}\put(21,
20){$g_5$}
\end{picture}

 The vertices   $g_0, g_1, g_3, g_5$ will be called odd and marked with    $\bullet$ on the graph,
 the vertices $g_2, g_4, g_5, g_6$  are even and indicated with $\circ$.
 The parameters of the  corresponding algebra
 $\mathcal{A}_{\alpha,\beta,\delta,\gamma}$ are enumerated according to the following picture:
$$
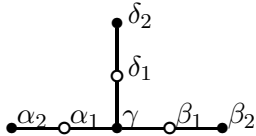
\begin{picture}(62,25)(-1,-1)
\thicklines
 \drawline(19,0)(19,9)
 \multiputlist(0,0)(10,0){\circle*{2},\circle{2},\circle*{2},\circle{2},\circle*{2}}
 \multiputlist(20,10)(0,10){\circle{2},\circle*{2}}
 \drawline(0,0)(8,0)
\drawline(10,0)(19,0)
\drawline(20,0)(28,0)\drawline(30,0)(38,0)\drawline(19,11)(19,19)
\put(0,1){$\alpha_2$} \put(10,1){$\alpha_1$} \put(20,1){$\gamma$}
\put(30,1){$\beta_1$} \put(40,1){$\beta_2$}
\put(21,10){$\delta_1$} \put(21,20){$\delta_2$}
\end{picture}
$$

We will write dimension and parameter vectors as
$(v_1,v_2,v_3,v_4,v_6,v_0)$. The minimal imaginary root for
$\widetilde{E_6}$ is the following vector $(1,2,1,2,1,2,3)$

Now we describe the root system for the diagram $\widetilde{E_6}$.
It consists of 72 series $\{\alpha_i +k \delta
\}_{k\in\mathbb{Z}}$ where $\alpha_i \in \Delta_f \cup -\Delta_f$

$\Delta_f = \{ (0, 0, 0, 0, 0, 0, 1), (0, 0, 0, 0, 1, 0, 0), (0,
0, 0, 0, 0, 1, 0), (0, 0, \ 0, 0, 0, 1, 1), \\ (0, 0, 0, 0, 1, 1,
0), (0, 0, 0, 0, 1, 1, 1), (0, 0, 1, 0, 0, \ 0, 0), (0, 0, 0, 1,
0, 0, 0), (0, 0, 0, 1, 0, 0, 1), \\ (0, 0, 0, 1, 0, 1, 1), \ (0,
0, 0, 1, 1, 1, 1), (0, 0, 1, 1, 0, 0, 0), (0, 0, 1, 1, 0, 0, 1),
(0, 0, \ 1, 1, 0, 1, 1), \\ (0, 0, 1, 1, 1, 1, 1), (0, 1, 0, 0, 0,
0, 0), (0, 1, 0, 0, 0, \ 0, 1), (0, 1, 0, 0, 0, 1, 1), (0, 1, 0,
0, 1, 1, 1), \\ (0, 1, 0, 1, 0, 0, 1), \ (0, 1, 0, 1, 0, 1, 1),
(0, 1, 0,
1, 0, 1, 2), (0, 1, 0, 1, 1, 1, 1), (0, 1, \ 0, 1, 1, 1, 2), \\
(0, 1, 0, 1, 1, 2, 2), (0, 1, 1, 1, 0, 0, 1), (0, 1, 1, 1, 0, \ 1,
1), (0, 1, 1, 1, 0, 1, 2), (0, 1, 1, 1, 1, 1, 1), \\ (0, 1, 1, 1,
1, 1, 2), \ (0, 1, 1, 1, 1, 2, 2), (0, 1, 1, 2, 0, 1, 2), (0, 1,
1, 2, 1, 1, 2), \\ (0, 1, \ 1, 2, 1, 2, 2), (0, 1, 1, 2, 1, 2, 3),
(0, 2, 1, 2, 1, 2, 3) \} $

The root system decomposes onto 7 C-series. Taking into account
the obvious symmetry of the graph $\widetilde{E_6}$ we need to
consider only three of them containing vectors $(1,0,0,0,0,0,0),
(0,1,0,0,0,0,0), (0,0,0,0,0,0,1)$ correspondingly. These C-series
 are

 $
 K_1=
 \{ (0, -2, -1, -2, -1, -2, -3), (0, -1, -1, -2, -1, -2, -3),
(0, -1, -1, -1, -1, \ -1, -1), \\ (0, -1, 0, 0, 0, 0, -1), (0, 0,
-1, -1, -1, -1, -1), (0, 0, 0, -1,  0, -1, -1), \\ (0, 0, 0, 1, 0,
1, 1), (0, 0, 1, 1, 1, 1, 1), (0, 1, 0, 0, 0, 0, \ 1), \\  (0, 1,
1, 1, 1, 1, 1), (0, 1, 1, 2, 1, 2, 3), (0, 2, 1, 2, 1, 2, 3) \}
+\delta \mathbb{Z}$

$ K_2 = \{  (0, -1, -1, -2, -1, -2, -2), (0, -1, -1, -1, -1, -1,
-2), (0, -1, 0, 0, 0, 0, \ 0), \\ (0, 1, 0, 0, 0, 0, 0), (0, 1, 1,
1, 1, 1, 2), (0, 1, 1, 2, 1, 2, 2)   \} +\delta \mathbb{Z}$

$ K_3= \{ (0, -1, 0, -1, 0, -1, -1), (0, 0, 0, 0, 0, 0, -1), (0,
0, 0, 0, 0, 0, 1), (0, \ 1, 0, 1, 0, 1, 1)  \} +\delta \mathbb{Z}$

\section{ Representations of algebras generated by
projections.}\label{repconnec}

Let us consider a tree $G$ with vertices  $\{ g_i,i=0,\ldots,
k+l+m\}$ and edges $\gamma_{g_i g_j}$,  see the figure.

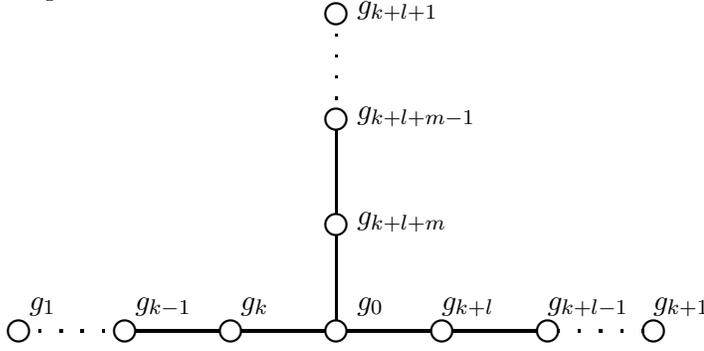
\setlength{\unitlength}{4pt}
\begin{picture}(70,32)(-1,-1)
\thicklines
 \drawline(29,1)(29,9)
\drawline(29,11)(29,19)
 \multiputlist(0,0)(10,0){\circle{2},\circle{2},\circle{2},\circle{2},\circle{2},\circle{2},\circle{2}}
 \multiputlist(30,10)(0,10){\circle{2},\circle{2},\circle{2}}
 \drawline(10,0)(18,0)
\drawline(20,0)(28,0) \drawline(30,0)(38,0)\drawline(40,0)(48,0)
\dottedline{2}(1,0)(7,0) \dottedline{2}(50,0)(58,0)
\dottedline{2}(29,22)(29,28)
 \put(31,10){$g_{k+l+m}$}\put(31,20){$g_{k+l+m-1}$}\put(31,30){$g_{k+l+1}$}
 \put(39,2){$g_{k+l}$}\put(49,2){$g_{k+l-1}$}\put(59,2){$g_{k+1}$}
\put(0,2){$g_{1}$}\put(10,2){$g_{k-1}$}\put(20,2){$g_{k}$}\put(31,2){$g_{0}$}
\end{picture}
\vspace{10pt}

We will establish a connection  between  $*$-representations of
the $*$-algebra  $\mathcal{A}_{G, \chi}$  and locally-scalar
representations of the graph  $G$, see~\cite{KPS}.

 Let us remark
that we can assume that each set  $M_1, M_2, \ldots, M_n$ contains
zero (this can be achieved by a translation). Henceforth,
$M_j=\alpha^{(j)}\cup\{0\}$. For three operators we will use
$\alpha$, $\beta$, $\delta$ instead of $\alpha^{(1)}$,
$\alpha^{(2)}$, $\alpha^{(3)}$. Thus vector $\chi$ will be
$(\alpha_1,\alpha_2,\ldots, \alpha_k, \beta_1, \beta_2,\ldots,
\beta_l, \delta_1,\delta_2,\ldots, \delta_m, \gamma)$.

\begin{definition}
An irreducible finite dimensional   *-representation  $\pi$ of the
algebra  $ \mathcal{A}_{G, \chi}$ such that
 $\pi(p_i)\not= 0 (1\le i\le k),\ \pi(q_j)
 \not= 0 (1\le j\le l),
 \ \pi(s_d)\not= 0 (1\le d\le m)$, and
 $\sum_{i=1}^k \pi(p_i)\not= I, \sum_{j=1}^l \pi(q_j)\not= I,
 \sum_{d=1}^m \pi(s_d)\not= I$,  will be called
  {\it non-degenerate}. By
 $\overline{\rep} \mathcal{A}_{G, \chi}$
 we will denote the full subcategory of non-degenerate
representations in the category ${\rep} \mathcal{A}_{G, \chi}$ of
$*$-representations of the $*$-algebra $\mathcal{A}_{G, \chi}$ in
the category $\mathcal{H}$ of Hilbert spaces.
\end{definition}

Let $\pi$ be a $*$-representation of $\mathcal{A}_{G, \chi}$
 on a Hilbert space  $H_0$. Put
  $P_i=\pi(p_i)$, $1\le i \le k$,
 $Q_j=\pi(q_j)$, $1\le j \le l$, $S_t=\pi(s_t)$, $1\le t \le m$.
 Let   $H_{p_i}= \image P_i$, $H_{q_j}= \image Q_j$, $H_{s_t}=
 \image S_t$. Denote by  $\Gamma_{p_i}$, $\Gamma_{q_j}$,
 $\Gamma_{s_t}$ the corresponding  natural isometries
 $H_{p_i}\to H_0$, $H_{q_j}\to H_0$,  $H_{s_t}\to H_0$.
 Then, in particular,  $\Gamma_{p_i}^*\Gamma_{p_i}
 =I_{H_{p_i}}$ is the  identity operator on  $H_{p_i}$ and
 $\Gamma_{p_i}\Gamma_{p_i}^*
 =P_i$. Similar equalities hold for the operators
 $\Gamma_{q_i}$ and $\Gamma_{s_i}$. Using  $\pi$  we construct
 a locally-scalar
representation   $\Pi$ of the graph   $G$.

  Let   $\Gamma_{ij}:H_j\to H_i$ denote the operator adjoint to   $\Gamma_{ji}:H_i\to H_j$, i.e. $\Gamma_{ij}=\Gamma_{ji}^*$. Put
\begin{gather*}
\Pi(g_0)= H^{g_0}=H_0,\\
 \Pi(g_k)= H^{g_k}= H_{p_1}\oplus H_{p_2}\oplus \ldots \oplus
H_{p_k},\\
\Pi(g_{k-1})= H^{g_{k-1}}=  H_{p_2}\oplus \ldots \oplus
H_{p_{k-1}}\oplus H_{p_k},\\
\Pi(g_{k-2})= H^{g_{k-2}}= H_{p_2}\oplus H_{p_3}\oplus \ldots
\oplus
H_{p_{k-1}},\\
\ldots.
 \end{gather*}

In these equalities the summands are omitted  from the left and
the right in turns. Analogously, we define subspaces  $\Pi(g_i)$
for $i=k+1,\ldots,k+l$ and  $i=k+l+1,\ldots,k+l+m$. Define  the
 operators  $\Gamma_{g_0,g_i}:H^{g_i}\to H^{g_0}$, where
$i\in\{k,k+l,k+l+m\}$, by the block-diagonal matrices
\begin{gather*}
\Gamma_{g_0,g_k}=\big[ \sqrt{\alpha_1} \Gamma_{p_1}|
\sqrt{\alpha_2} \Gamma_{p_2}|\ldots | \sqrt{\alpha_k} \Gamma_{p_k}
\big],\\
\Gamma_{g_0,g_{k+l}}=\big[ \sqrt{\beta_1} \Gamma_{q_1}|
\sqrt{\beta_2} \Gamma_{q_2}|\ldots | \sqrt{\beta_l} \Gamma_{q_l}
\big],\\
\Gamma_{g_0,g_{k+l+m}}=\big[ \sqrt{\delta_1} \Gamma_{s_1}|
\sqrt{\delta_2} \Gamma_{s_2}|\ldots | \sqrt{\delta_m} \Gamma_{s_m}
\big].\\
\end{gather*}
Now we define the representation  $\Pi$ on the edges
$\gamma_{g_0,g_k}$, $\gamma_{g_0,g_{k+l}}$,
$\gamma_{g_0,g_{k+l+m}}$ by the rule
\begin{gather*}
\Pi(\gamma_{g_0,g_k})=\{\Gamma_{g_0,g_k},\Gamma_{g_k,g_0} \},\\
\Pi(\gamma_{g_0,g_{k+l}})=\{\Gamma_{g_0,g_{k+l}},\Gamma_{g_{k+l},g_0} \},\\
\Pi(\gamma_{g_0,g_{k+l+m}})=\{\Gamma_{g_0,g_{k+l+m}},\Gamma_{g_{k+l+m},g_0} \}.\\
\end{gather*}

It is easy to see that
\[\Gamma_{g_0,g_k}
\Gamma_{g_k,g_0}+\Gamma_{g_0,g_{k+l}}
\Gamma_{g_{k+l},g_0}+\Gamma_{g_0,g_{k+l+m}}
\Gamma_{g_{k+l+m},g_0}=\gamma I_{H^{g_0}}.
\]

Let  $\mathcal{O}_{H,0}$  denote the operators from the zero space
to  $H$, and  $\mathcal{O}_{0,H}$  denote the  zero operator from
$H$ into the zero subspace. For the  operators $\Gamma_{g_j,g_i}
:H^{g_i}\to H^{g_j}$ with $i,j\not=0$, put
\begin{multline}\label{zved}
\Gamma_{g_{k-1},g_k}=\mathcal{O}_{0,H_{p_1}}\oplus
\sqrt{\alpha_{1} - \alpha_2} I_{H_{p_2}} \oplus \sqrt{\alpha_{1} -
\alpha_3} I_{H_{p_3}}\oplus \ldots \oplus \sqrt{\alpha_{1} -
\alpha_{k}}
I_{H_{p_{k}}}, \\
 \Gamma_{g_{k-1},g_{k-2}}=\sqrt{\alpha_{2}
- \alpha_k} I_{H_{p_2}} \oplus \sqrt{\alpha_{3} - \alpha_k}
I_{H_{p_3}} \oplus \ldots \oplus \sqrt{\alpha_{k-1} - \alpha_{k}}
I_{H_{p_{k-1}}}\oplus\mathcal{O}_{H_{p_k},0},\\
\Gamma_{g_{k-3},g_{k-2}}=\mathcal{O}_{0,H_{p_{2}}}\oplus
\sqrt{\alpha_{2} - \alpha_3} I_{H_{p_3}} \oplus \sqrt{\alpha_{2} -
\alpha_4} I_{H_{p_4}}\oplus \ldots \oplus \sqrt{\alpha_{2} -
\alpha_{k-1}}
I_{H_{p_{k-1}}},\\
\ldots\ldots\ldots\ldots\ldots\ldots\ldots\ldots\ldots\ldots.
\end{multline}

The corresponding operators for the rest of the edges of
 $G$ can be constructed analogously. One can
check that the  operators  $\Gamma_{g_{i},g_{j}}$, where
$\Gamma_{g_{i},g_{j}}=\Gamma_{g_{i},g_{j}}^*$, define a
locally-scalar representation of the graph  $G$ with the following
character  $f$:

\begin{gather}\label{char}
\begin{alignat*}{3}
 f(g_k) &= \alpha_1, &   \qquad f(g_{k+l}) &=
\beta_1, &  \qquad f(g_{k+l+m}) &= \delta_1, \\
 f(g_{k-1}) &=
\alpha_{1}-\alpha_{k}, & \qquad f(g_{k+l-1}) &=
\beta_{1}-\beta_{l}, &
\qquad f(g_{k+l+m-1})&= \delta_{1}-\delta_{m}, \\
 f(g_{k-2}) &= \alpha_{2}-\alpha_{k}, & \qquad f(g_{k+l-2}) &=
\beta_{2}-\beta_{l}, & \qquad f(g_{k+l+m-2}) &=
\delta_{2}-\delta_{m}, \\
   f(g_{k-3}) &= \alpha_{2}-\alpha_{k-1}, & \qquad f(g_{k+l-3}) &=
\beta_{2}-\beta_{l-1}, & \qquad f(g_{k+l+m-3}) &=
\delta_{2}-\delta_{m-1}, \\
 f(g_{k-4}) &= \alpha_{3}-\alpha_{k-1}, & \qquad f(g_{k+l-4}) &=
\beta_{3}-\beta_{l-1}, & \qquad f(g_{k+l+m-4}) &=
\delta_{3}-\delta_{m-1},\\
  \ldots &   & \ldots  & &  \ldots \\
\end{alignat*}
\end{gather}

$f(g_0)=\gamma$.

And vice versa, if a locally-scalar representation of the graph
$G$ with the  character $f(g_i)=x_i\in\mathbb{R}^*$ corresponds to
a non-degenerate representation of $\mathcal{A}_{G, \chi}$, then
one can check that
\begin{gather*}
\alpha_1=x_k,\\
\alpha_k=x_k-x_{k-1},
\alpha_{2}=x_k-x_{k-1}+x_{k-2},\\
\alpha_{k-1}=x_k-x_{k-1}+x_{k-2}-x_{k-3},\\
\alpha_{3}=x_k-x_{k-1}+x_{k-2}-x_{k-3}+x_{k-4},\\
\ldots.
\end{gather*}
Here $x_j=0$ if $j\le 0$. Analogously one can find   $\beta_j$ and
$\delta_t$. We will denote   $\Pi$ by $\Phi(\pi)$.

 Let  $\pi$  and $\widetilde{\pi}$ be non-degenerate representations of the algebra $\mathcal{P}_{\alpha,\beta,\delta,\gamma}$ and $C_0$ an intertwining operator for these representations; this is a morphism from  $\pi$ to  $\widetilde{\pi}$ in the category  $\rep G$),
 $C_0: H_0\to \widetilde{H}_0$,  $C_0 \pi= \widetilde{\pi} C_0$.  Put
\begin{gather*}
C_{p_i}= \widetilde{\Gamma}_{p_i}^* C_0 \Gamma_{p_i}, C_{p_i}:
H_{p_i}\to \widetilde{H}_{p_i},\   1\le i \le k,\\
C_{q_j}= \widetilde{\Gamma}_{q_j}^* C_0 \Gamma_{q_j}, C_{q_j}:
H_{q_j}\to \widetilde{H}_{q_j},\   k+1\le j \le k+l,\\
C_{s_t}= \widetilde{\Gamma}_{s_t}^* C_0 \Gamma_{s_t}, C_{s_t}:
H_{s_t}\to \widetilde{H}_{s_t}.\   k+l+1\le t \le k+l+m,\\
\ldots
\end{gather*}
 Put
\begin{gather*}
C^{(g_0)}=C_0: H^{(g_0)}\to \widetilde{H}^{(g_0)},\\
 C^{(g_k)}= C_{p_1}\oplus C_{p_2} \oplus \ldots  \oplus C_{p_k} :
H^{(g_k)}\to \widetilde{H}^{(g_k)},\\
 C^{(g_{k-1})}=  C_{p_2} \oplus \ldots  \oplus C_{p_{k-1}} \oplus C_{p_k}:
H^{(g_{k-1})}\to \widetilde{H}^{(g_{k-1})},\\
 C^{(g_{k-2})}= C_{p_2}\oplus C_{p_3} \oplus \ldots  \oplus C_{p_{k-1}} :
H^{(g_{k-2})}\to \widetilde{H}^{(g_{k-2})},\\
\ldots
\end{gather*}
 Analogously one can construct the operators  $C^{(g_i)}$ for  $i\in \{ k+l,\ldots, k+l+m
 \}$.  It is routine to check that the operators  $\{C^{(g_i)}\}_{0\le i\le
 k+l+m}$ intertwine the representations   $\Pi=\Phi(\pi)$  and
 $\widetilde{\Pi}=\Phi(\widetilde{\pi})$.  Put    $\Phi(C_0)=\{ C^{(g_i)}_{0\le i\le
 k+l+m}\}$.  Thus we have defined a functor   $\Phi:
 \overline{\rep} \mathcal{A}_{G, \chi}\to
 {\rep} G$, see~\cite{KPS}.   Moreover, the functor   $\Phi$  is univalent and full.
Let  $\widetilde{\rep}(G,d,f)$ be the full subcategory of
irreducible representations of  ${\rep}(G,d,f)$.
   $\Pi\in \ob \widetilde{\rep}(G,d,f)$, $f(g_i)= x_i\in \mathbb{R}^+,
 d(g_i)=d_i\in \mathbb{N}_0$, where $f$ is  the character of   $\Pi$, $d$ its dimension.
 It easy to verify that the representation $\Pi$ is isomorphic  (unitary equivalent) to an irreducible representation
 from the image of the functor $\Phi$ if and only if
\begin{gather}\label{nondeg}
1.\    0<x_1<x_2<\ldots < x_k; 0<x_{k+1}<x_{k+2}<\ldots <
x_{k+l};\\
  0<x_{k+l+1}<x_{k+l+2}<\ldots < x_{k+l+m}; \\
2.\    0< d_1<d_2<\ldots < d_k< d_0;  0< d_{k+1}<d_{k+2}<\ldots < \label{ineq1}\\
d_{k+l}< d_0;  0< d_{k+l+1}<d_{k+l+2}<\ldots < d_{k+l+m}<
d_0.\label{ineq2}
\end{gather}
(All matrices of the representation of the graph $G$,  except for
$\Gamma_{g_0,g_k}, \Gamma_{g_k,g_0}, \Gamma_{g_0,g_{k+l}}$, $
\Gamma_{g_{k+l}, g_0}$, $\Gamma_{g_0,g_{k+l+m}},
\Gamma_{g_{k+l+m}, g_0}$, can be brought    to the "canonical"
form (\ref{zved}) by admissible transformations. Then the rest of
the matrices will naturally be  partitioned into  blocks, which
gives the matrices $\Gamma_{p_i}, \Gamma_{q_i}, \Gamma_{s_i}$, and
thus the projections $P_i, Q_i, R_i$).  An irreducible
representation  $\Pi$ of the graph $G$ satisfying conditions
(\ref{nondeg})--(\ref{ineq2})  will be called {\it
non-degenerate}.
 Let
\begin{gather*}
\dim H_{p_i}= n_i, 1\le i\le k; \\
\dim H_{q_j}= n_{k+j}, 1\le j\le l; \\
\dim H_{s_t}= n_{k+l+t}, 1\le t\le m; \\
\dim H_{0}= n_{0}.
\end{gather*}
The vector ${n}= (n_0,n_1,\ldots , n_{k+l+m})$ is called the {\it
generalized dimension} of the representation  $\pi$ of the algebra
$\mathcal{A}_{G, \chi}$.  Let $\Pi=\Phi(\pi)$ for a non-degenerate
representation of the algebra $\mathcal{A}_{G, \chi}$,
${d}=(d_0,d_1,\ldots, d_{k+l+m})$ be the dimension of  $\Pi$.  It
is easy to see that
\begin{gather*}
n_1+n_2+\ldots + n_k = d_k,\\
n_2+\ldots + n_{k-1}+n_k = d_{k-1},\\
    n_2+\ldots + n_{k-1} = d_{k-2},\\
    n_3+\ldots + n_{k-2}+n_{k-1} = d_{k-3},\\
    \ldots
\end{gather*}
 Thus
\begin{multline}\label{dim}\\
n_1 =d_k-d_{k-1},\\
n_k =d_{k-1}-d_{k-2},\\
n_{2} =d_{k-2}-d_{k-3},\\
    \ldots \\
\end{multline}
  Analogously one can find  $n_{k+1}, \ldots , n_{k+l}$ from
  $d_{k+1}, \ldots, d_{k+l}$  and  $n_{k+l+1}, \ldots ,n_{k+l+m}$ from  $d_{k+l+1}, \ldots, d_{k+l+m}$

 Denote by  $\overline{\rep}G$  the full subcategory in   $\rep
 G$ of non-degenerate locally-scalar representations of the graph   $G$.
As a corollary of the previous arguments we obtain the following
theorem.
\begin{theorem}
Let  $\mathcal{A}_{G, \chi}$ be associated with a graph  $G$. The
functor  $\Phi$ is an equivalence of the categories
$\overline{\rep}{ \mathcal{A}_{G, \chi}}$
 of  non-degenerate *-representations of the algebra
 $\mathcal{A}_{G, \chi}$ and the
category $\overline{\rep}G$ of non-degenerate locally-scalar
representations of the graph  $G$.
\end{theorem}

Let us define the Coxeter functors for the  $*$-algebras
$\mathcal{A}_{G, \chi}$, by putting $\overset{\circ}{\Psi}=
\Phi^{-1}\overset{\circ}{F} \Phi $ and $\overset{\bullet}{\Psi}=
\Phi^{-1}\overset{\bullet}{F} \Phi$. Now we can use the results
of~\cite{roiter} to give a description of representation of the
*-algebra $\mathcal{A}_{G, \chi}$.

Note that  to find  formulae of the lacaly-scalr representations
of a given extended Dynkin graph we need to consider two
principally different cases: the case when the vector of
generalized dimension is a real root and the case when it is a
imaginary root. In the letter case the vector of parameters $\chi$
must satisfy (in order for the representations to exist) a certain
linear relation ("traces equality"). Hence $\chi$ must belong to a
ceratin hyperplane $h_G$ which depends only on the graph $G$. A
simple calculation yields that for extended Dynkin graphs
$\widetilde{D}_4$, $\widetilde{E}_6$, $\widetilde{E}_7$,
$\widetilde{E}_8$ these hyperplanes are the following:

\begin{tabular}{|p{2cm}|p{2cm}|p{2cm}|p{2cm}|}
  \hline
  $\widetilde{D}_4$ & $\widetilde{E}_6$ & $\widetilde{E}_7$ & $\widetilde{E}_8$ \\
  \hline
  $\alpha_1+\beta_1+\delta_1 +
\eta_1 = 2 \gamma$ & $\alpha_1+\alpha_2+\beta_1+\beta_2
+\delta_1+\delta_2 = 3 \gamma$ & $\alpha_1+\alpha_2+\alpha_3+
\beta_1+\beta_2 +\beta_3+ 2\delta_1 = 4 \gamma$
& $2 (\alpha_1+\alpha_2)+ \beta_1+\beta_2 +  \beta_3+ \beta_4 +3\delta_1 = 6 \gamma$  \\
  \hline
\end{tabular}

It is know (see~\cite{melit}) that in case $\chi \in h_G$ the
dimension of any irreducible representation is bounded (by 2 for
$\widetilde{D}_4$, by 3 for $\widetilde{E}_6$, by 4 for
$\widetilde{E}_7$ and by 6 for $\widetilde{E}_8$). Thus in case of
the hyperplane we can describe the set of admissible parameters
$\chi$ using Horn's inequalities. In case $\chi\not\in h_g$
 the dimension of any  irreducible locally-scalar representation
 is a real root. In what follows we will relay on the following result due to
V.Ostrovskij \cite{ostrovskij}
\begin{theorem}
Let $\pi$ be an irreducible *-representation of the algebra
$\mathcal{A}_{G,\chi}$ associated with Extended Dynkin graph $G$
and $\widehat{\pi}$ corresponding representation of the graph $G$.
Then either generalized dimension $d$ of $\widehat{\pi}$ is a
singular root or vector-parameter $\chi\in h_G$.
\end{theorem}
Hence we will solve the spectral problem if we describe the
parameters $\chi$ for which there are locally-scalar
representations with vector of generalized dimension being real
singular roots and parameters belonging to  hyperplane $h_G$ for
which there exist representations of the algebra. Firstly we will
consider the case $\chi\not\in h_G$.

In the next section we will do the following: we know how to
construct all irreducible locally-scalar representations of Dynkin
graphs with the aid of Coxeter reflection functors starting from
the simplest ones. In particular, we can find their dimensions and
characters \cite{roiter}. Next we single out non-generate
representations and apply the  equivalence functor  $\Phi$, see
Theorem~1.
\section{Algebras associated with extended Dynkin graphs.}\label{s2}

Thus in "generic" situation the powers of the Coxeter map
$C^k=(\overset{\bullet}{C}\overset{\circ}{C})^k$ acts on vectors
(in quiver notations) by the following formulas:

if $k\equiv 0(\mod \ 3)$
\[\frac{1}{12}
\left(%
\begin{array}{ccccccc}
  11+3 (-1)^k & 4 & -1+3 (-1)^k & 4 & -1+3 (-1)^k & 4 & 3 (3-(-1)^k) \\
    -4 & 4 & -4 & -8 & -4 & -8 & -12   \\
 -1+3 (-1)^k & 4 & 11+3 (-1)^k & 4 & -1+3 (-1)^k & 4 & 3 (3-(-1)^k) \\
   -4 & -8 & -4 & 4 & -4 & -8 & -12   \\
 -1+3 (-1)^k & 4 & -1+3 (-1)^k & 4 & 11+3 (-1)^k & 4 & 3 (3-(-1)^k)  \\
   -4 & -8 & -4 & -8 & -4 & 4 & -12   \\
 3 (3-(-1)^k) & 12 & 3 (3-(-1)^k) & 12 & 3 (3-(-1)^k) & 12 & 3 (9+(-1)^k)
\end{array}%
\right)
\]

if $k\equiv 1(\mod \ 3)$
\[\frac{1}{4}
\left(%
\begin{array}{ccccccc}
1+(-1)^k& 4 &1+(-1)^k& 0& 1+(-1)^k& 0& 3-(-1)^k \\
 -4&-4&0&0&0&0&-4 \\
  1+(-1)^k& 0& 1+(-1)^k& 4& 1+(-1)^k& 0& 3-(-1)^k \\
 0& 0& -4& -4& 0& 0& -4 \\
  1+(-1)^k& 0& 1+(-1)^k& 0& 1+(-1)^k& 4& 3-(-1)^k \\
 0& 0& 0& 0& -4& -4& -4 \\
  3-(-1)^k&4&3-(-1)^k&4&3-(-1)^k&4&9+(-1)^k
\end{array}%
\right)
\]

if $k\equiv 2(\mod \ 3)$
\[ \frac{1}{12}
\left(%
\begin{array}{ccccccc}
-5+3 (-1)^k & -4 & 7+3 (-1)^k & 8 & 7+3 (-1)^k & 8 & 3 (3-(-1)^k)  \\
  4 & -4 & -8 & -4 & -8 & -4 & -12   \\
 7+3 (-1)^k & 8 & -5+3 (-1)^k & -4 & 7+3 (-1)^k & 8 & 3 (3-(-1)^k)  \\
   -8 & -4 & 4 & -4 & -8 & -4 & -12   \\
 7+3 (-1)^k & 8 & 7+3 (-1)^k & 8 & -5+3 (-1)^k & -4 & 3 (3-(-1)^k)  \\
   -8 & -4 & -8 & -4 & 4 & -4 & -12   \\
 3 (3-(-1)^k) & 12 & 3 (3-(-1)^k) & 12 & 3 (3-(-1)^k) & 12 & 3
(9+(-1)^k)
\end{array}%
\right)
\]

Let us consider the first C-series which contains simple root
$(1,0,0,0,0,0,0)$. It decomposes onto 12 $\delta$-series. The map
$C^6$ takes simple representation of the quiver to the
non-degenerate one. To obtain further representations of this
series we need to apply $C^k$ written explicitly above.

Let $M_d$ denote transition matrix from generalized dimension of
the quiver to the generalized  dimension of the algebra and $M_f$
denote the transition matrix from algebra parameters to quiver
parameters.

We will write $v\ge 0$ meaning that $v_j > 0$ for $1\le j\le 6$
and $v_7 = 0$.

\begin{theorem}\label{ineqe6}

Consider the following vectors $v_1=(1, 0, 0, 0, 0, 0, 0), v_2=(1,
1, 0, 0, 0, 0, 0), v_3=(0, 1, 0, 0, 0, 0, 1), v_4=(0, 0,  0, 1, 0,
1, 1), v_5=(0, 0, 1, 1, 1, 1, 1), v_6=(0, 1, 1, 1, 1, 1, 1),
v_7=(1, 1, 0, 1, 0, \ 1, 2), v_8=(1, 2, 0, 1, 0, 1, 2), v_9=(1, 2,
1, 1, 1, 1, 2), v_{10}=(1, 1, 1, 2, 1, 2, 2), v_{11}=(0, 1, 1, 2,
1, 2, 3), v_{12}=(0, 2, 1, 2, 1, 2, 3)$. Put $d_k= v_{k \mod\ 12}+
[\frac{k}{12}]\delta$.

 The algebra  $\mathcal{A}_{\widetilde{E_6},\chi}$  associated with the  Dynkin
graph  $\widetilde{E_6}$ has an irreducible non-degenerate
representation  in generalized dimensions $M_d d_k$ for all $k\ge
15$ iff the parameters $\chi$ satisfy the following condition:
$D_k \chi \ge 0$ where matrix
$$D_k=D_1\begin{cases}  C^{6-s} M_f  &
\text{if $k=2s$}\\
C^{6-s}\overset{\circ}{C} M_f  & \text{if $k=2s+1$}
\end{cases}
$$ and $$D_1 = \left(%
\begin{array}{ccccccc}
  3 & -3 & 1 & -3 & 1 & -3 & 5 \\
  1 & -1 & 1 & -2 & 1 & -1 & 2 \\
  3 & -3 & 2 & -3 & 1 & -2 & 4 \\
  1 & -1 & 1 & -1 & 1 & -2 & 2 \\
  3 & -3 & 1 & -2 & 2 & -3 & 4 \\
  5 & -4 & 2 & -4 & 2 & -4 & 6 \\
  2 & -2 & 1 & -2 & 1 & -2 & 3 \\
\end{array}%
\right)$$
\end{theorem}

\begin{theorem}\label{ineqe6_2}

Consider the following vectors $v_1=(0, 1, 0, 0, 0, 0, 0), v_2=(1,
1, 0, 0, 0, 0, 1), v_3=(1, 1, 0, 1, 0, 1, 1), v_4=(0, 1, 1, 1, 1,
1, 2), v_5=(0, 1, 1, 2, 1, 2, 2), v_6=(1, 1, 1, 2, 1, 2, 3)$. Put
$d_k= v_{k\mod 6}+ [\frac{k}{6}]\delta$.

 The algebra  $\mathcal{A}_{\widetilde{E_6},\chi}$  associated with the  Dynkin
graph  $\widetilde{E_6}$ has an irreducible non-degenerate
representation  in generalized dimensions $M_d d_k$ for all $k\ge
8$ iff the parameters $\chi$ satisfy the following condition: $D_k
\chi \ge 0$ where matrix
$$D_k=D_2\begin{cases}  C^{3-s} M_f  &
\text{if $k=2s$}\\
C^{3-s}\overset{\circ}{C} M_f  & \text{if $k=2s+1$}
\end{cases}
$$ where
$$ D_2= \left(%
\begin{array}{ccccccc}
  0 & 0 & 1 & -1 & 1 & -1 & 1 \\
  0 & 0 & 0 & -1 & 0 & 0 & 1 \\
  1 & -1 & 0 & -1 & 1 & -1 & 2 \\
  0 & 0 & 0 & 0 & 0 & -1 & 1 \\
  1 & -1 & 1 & -1 & 0 & -1 & 2 \\
  2 & -1 & 1 & -2 & 1 & -2 & 3 \\
  1 & -1 & 1 & -2 & 1 & -2 & 3 \\
\end{array}%
\right) $$

\end{theorem}

\begin{theorem}\label{ineqe6_3}

Consider the following vectors $v_1=(0, 0, 0, 0, 0, 0, 1), v_2=(0,
1, 0, 1, 0, 1, 1), v_3=(1, 1, 1, 1, 1, 1, 2), v_4=(1, 2, 1, 2, 1,
2, 2)$. Put $d_k= v_{k\mod 4}+ [\frac{k}{4}]\delta$.

 The algebra  $\mathcal{A}_{\widetilde{E_6},\chi}$  associated with the  Dynkin
graph  $\widetilde{E_6}$ has an irreducible non-degenerate
representation  in generalized dimensions $M_d d_k$ for all $k\ge
5$ iff the parameters $\chi$ satisfy the following condition: $D_k
\chi \ge 0$ where matrix
$$D_k=D_3\begin{cases}  C^{s-2} M_f  &
\text{if $k=2s$}\\
C^{s-2}\overset{\bullet}{C} M_f  & \text{if $k=2s+1$}
\end{cases}
$$ where
$$ D_3= \left(%
\begin{array}{ccccccc}
  -1 & 1 & 0 & 0 & 0 & 0 & 0 \\
  -1 & 1 & 0 & 1 & 0 & 1 & -1 \\
  0 & 0 & -1 & 1 & 0 & 0 & 0 \\
  0 & 1 & -1 & 1 & 0 & 1 & -1 \\
  0 & 0 & 0 & 0 & -1 & 1 & 0 \\
  0 & 1 & 0 & 1 & -1 & 1 & -1 \\
  -1 & 2 & -1 & 2 & -1 & 2 & -2 \\
\end{array}%
\right) $$

\end{theorem}

\begin{theorem}\label{E6_summ}
Up to symmetry all irreducible non-degenerate *-representations of
the algebra $\mathcal{A}_{\widetilde{E_6},\chi}$ with vector of
parameters not on a hyperplane are described in
theorems\ref{ineqe6}, \ref{ineqe6_2},\ref{ineqe6_3}.
\end{theorem}

\section{The hyperplane case. }

 Applying Horn's inequalities we get the
following theorem:
\begin{theorem}\label{E6_hyp}
The algebra  $\mathcal{A}_{\widetilde{E_6},\chi}$  associated with
the  Dynkin graph  $\widetilde{E_6}$ has *-representation  in
generalized dimensions $(1,1;1,1;1,1;3)$  iff the parameters
$\chi$ satisfy the following conditions:

$2 (\alpha_1 +\beta_1 )>\alpha_2 +\beta_2 + \delta_1  +\delta_2, 2
(\alpha_1  +\delta_1 )> \alpha_2  +\beta_1  +
\beta_2  +\delta_2,\\
2 (\beta_1  +\delta_1  )>\alpha_1  +\alpha_2  +\beta_2  +
\delta_2, \alpha_1  +\alpha_2  +\beta_1  +\delta_1 >2 (\beta_2
 +\delta_2),\\
2 (\alpha_2  +\beta_2  +\delta_1)>\alpha_1  +\beta_1
 +\delta_2,
\alpha_1  +\beta_1  +\beta_2  +\delta_1 >2 (\alpha_2
+\delta_2),\\
2 (\alpha_2  +\beta_1  +\delta_2  )>\alpha_1  +\beta_2  +
\delta_1, 2 (\alpha_1  +\beta_2  +\delta_2  )>\alpha_2
 +\beta_1  +\delta_1,\\
  \alpha_1  +\alpha_2  +\beta_1  +\beta_2  +\delta_2  >2 \delta_1,
  \alpha_1  +\beta_1  +\delta_1  +\delta_2  >2 (\alpha_2  +\beta_2),\\
  \alpha_1  +\alpha_2  +\beta_2  +\delta_1  +\delta_2  >2 \beta_1,
  \alpha_2  +\beta_1  +\beta_2  +\delta_1  +\delta_2  >2 \alpha_1.
  $
\end{theorem}
The following theorem gives a complete solution to the spectral
problem in case of Dynkin graph $\widetilde{E}_6$.
\begin{theorem} A non-degenerate *-representation of the algebra
$\mathcal{A}_{\widetilde{E_6},\chi}$ exists iff $\chi$ satisfies
at the conditions of at least one of the theorems \ref{ineqe6},
\ref{ineqe6_2}, \ref{ineqe6_3}, \ref{E6_hyp}.
\end{theorem}
\begin{proof}
The case $\chi \not\in h_{\widetilde{E_6}}$ is covered by theorems
\ref{ineqe6}-\ref{ineqe6_3}.

 If $\chi\in h_{\widetilde{E_6}}$ then for any irreducible
 *-representation $\pi$ of the algebra
 $\mathcal{A}_{\widetilde{E_6},\chi}$ the generalized dimension of
 the corresponding representation $\hat{\pi}$ of
 the graph $\widetilde{E_6}$ is the minimal imaginary
 root $\delta = (1,2;1,2;1,2;3)$ (then dimension of $\pi$ is
 $(1,1;1,1;1,1;3)$)  or is a real root $d$.

 We have already described the case of singular roots in theorems
\ref{ineqe6}-\ref{ineqe6_3}. But $d$ can not be regular. Because
otherwise application of the Coxeter functors  would produce
irreducible representations in all dimensions of the form $d+k
\delta$ where $k$ is a positive integer. This would be a
contradiction with the fact that
$\mathcal{A}_{\widetilde{E_6},\chi}$ is PI algebra for $\chi \in
h_{\widetilde{E_6}}$.
 Thus we have exhausted all the possibilities.

\end{proof}

\end{document}